\documentclass[12pt,leqno]{article}
\usepackage{latexsym,amsmath,amssymb,epsfig}
\title{Quotients of a Universal Locally Projective Polytope
of type $\{5,3,5\}$
{ }}
\author{Michael I.~Hartley\\ {\small KDU College}\\
{\small Jalan SS 22/41, Damansara Jaya}\\
{\small Selangor, 47400, Malaysia}\\
{\small michael@kdu.edu.my} \\ \\
Dimitri Leemans\\ {\small D\'epartement de Math\'ematiques}\\
{\small Universit\'e Libre de Bruxelles}\\
{\small C.P.216~-- G\'eom\'etrie}\\
{\small Boulevard du Triomphe}\\
{\small B-1050 Bruxelles}\\
{\small dleemans@ulb.ac.be}}
\date{}

\advance\voffset by-20mm
\makeatletter
\def\@eqnnum{\hbox to .01pt{}\rlap{\bf \hskip -
\displaywidth\theequation}}
\makeatother
\renewcommand{\theequation}{\thesection.\arabic{equation}}
\newcommand{\seq}{\setcounter{equation}{0}}
\newenvironment{proof}{\par \noindent {\em
Proof\/}:\quad}{\hspace*{\fill} $\Box$ \par \vspace*{1ex}}
\newcommand{\bgpf}{\begin{proof}}
\newcommand{\ndpf}{\end{proof}}
\newenvironment{reason}{\par \noindent {\em
Reason\/}:\quad}{\hspace*{\fill} $\Box$ \par \vspace*{1ex}}
\newcommand{\bgrs}{\begin{reason}}
\newcommand{\ndrs}{\end{reason}}

\newtheorem{theorem}[equation]{\hspace{-0.3em}}
\newcommand{\bgth}{\begin{theorem}}
\newcommand{\ndth}{\end{theorem}}

\newcommand{\thm}{{\sc Theorem}\quad}

\newcommand{\bgeq}{\begin{equation}}
\newcommand{\bgqy}{\begin{eqnarray*}}
\newcommand{\bgry}{\begin{array}}
\newcommand{\ndeq}{\end{equation}}
\newcommand{\ndqy}{\end{eqnarray*}}
\newcommand{\ndry}{\end{array}}
\def\Aut{\mathop{{\rm Aut}}}

\newcommand{\CK}{{\cal K}}
\newcommand{\CL}{{\cal L}}

\newcommand{\CP}{{\cal P}}
\newcommand{\CQ}{{\cal Q}}

\newcommand{\scl}[1]{\langle{#1}\rangle}
\begin{document}
\maketitle
\centerline{\it In Memory of H. S. M. ``Donald'' Coxeter, 1907--2003}
\begin{abstract}
This article examines the universal polytope $\CP$ (of type $\{5,3,5\}$) whose
facets are dodecahedra, and whose vertex figures are hemi-icosahedra.
The polytope is proven to be finite, and the structure of its group
is identified. This information is used to classifiy the quotients
of the polytope. A total of 145 quotients are found, including 69
section regular polytopes with the same facets and vertex figures as $\CP$.
\end{abstract}

%
%
%
%
\section{Introduction}
\label{intro}
\seq

This article may be seen in three different ways. On the one hand, it was
inspired by, and solves, a problem in the theory of polytopes. A certain
regular polytope is proven to exist and be finite, and its quotient
polytopes are discovered. Being about polytopes, it is also a paper
about geometries or buildings, since (for example)
every regular polytope is a thin regular geometry with a string diagram.
Finally, it may also be regarded as an article about Coxeter groups. A certain
quotient of a hyperbolic Coxeter group is proven to be finite, and its structure
as a group is identified.

The polytope studied is the universal polytope $\CP$
whose facets are dodecahedra
and whose vertex figures are hemi-icosahedra. This polytope is
locally projective, that is, it is section
regular with facets and vertex figures either spherical or projective, and
not both spherical. The polytope, of type $\{5,3,5\}$, covers every
polytope with facets and vertex figures of these types.
A classification of its quotients is therefore
important since it encompasses
a classification of the locally projective polytopes of type $\{5,3,5\}$~--
every such polytope or its dual is a quotient of $\CP$. The automorphism
group of this polytope is the quotient of the Coxeter group
$[5,3,5]=\scl{s_0,s_1,s_2,s_3}$ by the normal subgroup generated by
all quotients of $(s_1s_2s_3)^5$.

The results of the article is expressed mainly in the language
of abstract polytopes, since that is the setting which gave rise to the
results, and (in the authors' opinion) in which they seem the most
significant.
It is assumed from this point on that the reader is familiar with the basic concepts of
abstract polytopes and their quotients. The most important reference in the theory of
abstract regular polytopes is \cite{msbook}. An abstract polytope is a poset
satisfying certain properties that are satisfied by the face-lattices of
classical polytopes, as well as honeycombs in euclidian and hyperbolic space,
and other objects. The theory therefore encompasses all of the latter, as well
as a rich assortment of polytopes on other spaceforms and objects for which
it is difficult to define a natural topology.

The locally projective polytopes fall into the latter class. A (section) regular
polytope may be defined to be {\it locally X} if its minimal nonspherical sections have topological
type X (see \cite{mslocproj}). Alternatively but not equivalently, some authors
define a locally X polytope to be one whose facets and vertex figures
are either spherical or X, but not both spherical (see for example
\cite{mstwis}). Note that the former
definition subsumes the latter, and that in rank 4 the two definitions
are equivalent, since all polytopes of rank 2 or less are spherical.

A regular polytope is one whose automorphism group acts transitively on
its set of flags (that is, maximal totally ordered subsets). The most
important result in the study of abstract regular polytopes is that they
are in one to one correspondence with so-called string C-groups, that is
groups generated by involutions $s_0,s_1,\dots,s_{n-1}$ where first of
all $s_is_j=s_js_i$ whenever $i\neq j,j\pm1$, and secondly,
$\scl{s_i:i\in I}\cap\scl{s_i:i\in J}=\scl{s_i:i\in I\cap J}$ for any
$I, J\subseteq\{0,\dots,n-1\}$.
The automorphism groups of regular polytopes are always C-groups, and from
any C-group the corresponding polytope may be reconstructed.

Given a polytope $\CP$ and a subgroup $N$ of its automorphism group, the elements
of $\CP$ are partitioned into orbits $\{F\cdot N:F\in\CP\}$ by $N$.
We may define a poset $\CP/N$ on these orbits, the {\it quotient} of $\CP$ by $N$,
in a natural way, letting
$F\cdot N\leq G\cdot N$ if $F\leq G$ (that is, $F\cdot N\leq G\cdot N$ if and
only if there exists $F'\in F\cdot N$ and $G'\in G\cdot N$ with $F'\leq G'$).
When $\CP$ is regular, the conditions on $N$ for which $\CP/N$ is again a polytope
are well-known (see \cite{mcq}, or Section 2D of \cite{msbook}).

In fact, every polytope may be written as a quotient of some regular polytope
(see \cite{har1} and \cite{har2}), and the quotients of a regular polytope $\CP$ are in
one-to-one correspondence with conjugacy classes of so-called ``semisparse''
subgroups of $\Aut(\CP)$.

Let $\CK$ and $\CL$ be regular polytopes. If there exists a polytope $\CQ$
whose facets are of type $\CL$ and whose vertex figures are of type $\CK$,
then there is a {\it universal} such polytope $\CP$, denoted $\{\CL,\CK\}$
(with automorphism group $[\CL,\CK]$),
which ``covers'' all other such polytopes in the sense that they are quotients
of $\CP$. This was shown in Theorem 2.5 of \cite{har6}, and much earlier for the case
when $\CQ$ is regular in \cite{amalg} (see Theorem 4A2 of \cite{msbook}).
The search for polytopes with particular facets and vertex figures therefore
usually follows the following pattern.
First the universal such polytope is discovered (if not already known). Secondly,
the quotients of this universal polytope are sought. In \cite{har6}, various
results about semisparse subgroups were uncovered that facilitate this process.

Let $W=\scl{s_0,\dots,s_{n-1}}$ be the group of a regular polytope $\CP$,
let $H_{n-1}=\scl{s_0,\dots,s_{n-2}}$ be the group of its facets, and
let $H_0=\scl{s_1,\dots,s_{n-1}}$ be the group of its vertex figures.
The key result from \cite{har6} that we shall use here is as follows.
If the vertex figures of $\CP$ have no proper quotients,
then subgroups $N$ of $W$ are semisparse if and only $N^w\cap H_0H_{n-1}$
is semisparse in $H_{n-1}$ for all conjugates $N^w$ of $N$ in $W$ (see
Theorem 2.7 of \cite{har6}).

In this article, we are interested in polytopes whose facets are dodecahedra,
and whose vertex figures are hemi-icosahedra. If any such polytopes exist
at all, they are quotients of a universal such polytope with automorphism group
$W=\scl{s_0,s_1,s_2,s_3}$, satisfying $s_0^2=s_1^2=s_2^2=s_3^2=1$,
$(s_0s_1)^5=(s_1s_2)^3=(s_2s_3)^5=1$, $s_0s_2=s_2s_0$, $s_1s_3=s_3s_1$,
$s_0s_3=s_3s_0$ and $(s_1s_2s_3)^5=1$. For the remainder of this article,
$W$ shall denote this group. It is shown in Section~\ref{twoquos}
that $W$ is a C-group, the group of a well-defined polytope $\CP$ whose facets and
vertex figures are as desired. This is done by exhibiting an example of
another (smaller) such polytope, which must therefore be a proper quotient
of $\CP$.

In Section~\ref{struct}, it is shown that $W$ (and therefore $\CP$)
is in fact finite, and is the direct product of two large simple groups,
the Janko group $J_1$ and the projective special linear group $L_2(19)$.
In Section~\ref{allquos}, the remaining quotients are discovered and
tabulated.

Let $H_0=\scl{s_1,s_2,s_3}$ and $H_3=\scl{s_0,s_1,s_2}$.

\section{Two Quotients of $\CP$.}
\label{twoquos}
\seq

Coxeter, in \cite{cox57}, discovered a self-dual locally projective polytope $\CP''$
with 57 hemi-dodecahedral facets (it was also constructed in \cite{cruyce}).
We shall call this polytope the
{\it 57-cell}. In \cite{har4} it was shown that Coxeter's
57-cell has no proper quotients. Its group $W''$ is generated by $\{s''_0,s''_1,s''_2,s''_3\}$
with the relations of $[5,3,5]$ and the additional relations $(s''_0s''_1s''_2)^5=(s''_1s''_2s''_3)^5=1$.
Moreover $W''$ is isomorphic to the simple group $L_2(19)(=PSL_2(19))$, of order $3420=2^2.3^2.5.19$. It is
a quotient of $W=\scl{s_0,s_1,s_2,s_3}$ by the normal subgroup $N''$
generated by all conjugates of $(s_0s_1s_2)^5$.
The simple group $L_2(19)$ has a permutation presentation on 20 points.
Readers interested to see this permutation presentation should download
this article's ``auxiliary information'' available at \cite{535aux}. (See also the notes
at the end of Section~\ref{allquos}.)

For the remainder of this article,
let $\omega=(s_0s_1s_2)^5$, so that $N''=\scl{\scl{\omega}}$. The action of
$\omega$ (as an automorphism) on $\CP$ is to move a (base) flag $\Phi$
to the ``opposite'' flag of the facet contained in $\Phi$.

Define $\nu=\omega s_3$, and let $N'=\scl{\scl{\nu^3}}$ be the normal subgroup of $W$
generated by all conjugates of $\nu^3$. The group $W'\cong W/N'$ may be taken
to be the group $\scl{s'_0,s'_1,s'_2,s'_3}$ whose generators
satisfy all the relations of $[5,3,5]$ as well as the additional
relations $(s'_1s'_2s'_3)^5=((s'_0s'_1s'_2)^5s'_3)^3=1$. A computer
algebra package \cite{gap4} was used to analyse the group $W'$. It was
found to be isomorphic to the Janko group $J_1$, a sporadic finite
simple group of order $175560=2^3.3.5.7.11.19$. It was also checked
that $W'$ is a C-group, that $H'_3=\scl{s'_0,s'_1,s'_2}$ has order
120, that $H'_0=\scl{s'_1,s'_2,s'_3}$ has order 60, and therefore
that $W'$ is the group of a regular polytope $\CP'$ with dodecahedral
facets and hemi-icosahedral vertex figures.

The simplest permutation representation for $J_1$ is a permutation action on
266 points (see \cite{janko} and \cite{liv} for more information).
Again, readers interested in an example of this permutation
action are referred to the auxiliary information for this article. (See the
notes at the end of Section~\ref{allquos} for details.)

This polytope $\CP'$ has $1463$ dodecahedral facets,
and twice that number of vertices. It is instructive to
classify the quotients of $\CP'$. Theorem 2.7 of \cite{har6}
applies, so a subgroup $K\leq W'$ is semisparse if and only if
$K^w\cap H'_0H'_3=\{1\}$ or $\{1,(s'_0s'_1s'_2)^5\}$.
In $J_1$, all elements of order
$2$ are conjugate. It follows that if $K$ contains $(s'_0s'_1s'_2)^5$,
it has a conjugate containing $s'_0$, which would contradict the
semisparseness of $K$. Therefore, $\CP'/K$ cannot have any
hemi-dodecahedral facets~-- all its facets must be dodecahedra, and
$K^w\cap H'_0H'_3=\{1\}$ for all conjugates $K^w$ of $K$ in $W'$.

The Sylow $p$-subgroups of any group are all conjugate. For $p=3$, $5$, $7$, $11$ or
$19$, the Sylow $p$-subgroups of $J_1$ are cyclic of order $p$.
If $H'_0H'_3$ contains an element of order $p$, then $K$ can not, otherwise
there would exist a conjugate of $K$ intersecting $H'_0H'_3$ nontrivially.
It may be verified, however, that $H'_0H'_3$ does contain elements of each of
these orders: for example $s'_1s'_2$ has order $3$, $s'_0s'_1$ has order 5,
$s'_3s'_1s'_2s'_1s'_0$ has order 19, $s'_3s'_2s'_1s'_0s'_1s'_0$ order 11 and
$s'_1s'_3s'_2s'_1s'_3s'_0s'_1s'_0s'_2s'_1s'_0$ order 7.
Therefore the polytope $\CP'$, like $\CP''$, has no proper quotients.

\section{The Group Structure of $W$.}
\label{struct}
\seq

We can now investigate the structure of the group $W=\scl{s_0,s_1,s_2,s_3}$.
An important result about $W$ is the following.

\bgth
\label{WisC}
\thm $W$ is a C-group, and the polytope $\CP$ has dodecahedral facets
and hemi-icosahedral vertex figures.
\ndth

\bgpf
Proposition 4A8 of \cite{msbook}, combined with either of the examples of the previous
section, show that the universal polytope $\{\{5,3\},\{3,5\}_5\}$ exists,
and that the presentation of its automorphism group is just that of $W$.
\ndpf

We already know that $W$ has two normal subgroups $N'$ and $N''$,
of index $175560$ and $3420$ respectively. This information may be
used to discover the structure of $W$.

\bgth
\label{WisJL}
\thm $W\cong J_1\times L_2(19)$.
\ndth

\bgpf
$W$ has a subgroup
$L=\scl{s_2s_1,(s_0s_1)^2s_2s_3s_0,s_0s_1s_0s_3s_2s_3s_1}$ of
index 20. That the index is 20 may be verified
 using the Todd-Coxeter coset enumeration technique, either by
 hand or using a computer.

Let $x_1$, $x_2$ and $x_3$ be the three generators of $L$ in the order given,
and let $x=x_2^{-1}x_3x_1$ and $y=x_1x_3^{-1}$. Then $x$ and $y$
generate $L$, and satisfy
the relations $y^3=xyxy^{-1}xyx^{-1}yx=1$. This may be shown either
laboriously by hand, or by using a computer. Finally, it may be shown
that the group $\scl{x,y:y^3=xyxy^{-1}xyx^{-1}yx=1}$ is finite, of order
30020760. 
Since $L$ is a quotient of this group, it follows that $W$ has order at most
$20\times 30020760 = 3420\times 175560$. In particular, $W$ is finite.
On the other hand, since $W$ has quotients isomorphic to the simple groups $L_2(19)$ and $J_1$,
we deduce that $W$ has order at least $3420\times 175560$, so in fact
this number is exactly the order of $W$, and $L$ equals
$\scl{x,y:y^3=xyxy^{-1}xyx^{-1}yx=1}$.

The details of this derivation of the structure of $W$ were obtained
using GAP version 4 (see \cite{gap4}). In particular, the neat presentation
for $L$ was obtained in GAP using the Reduced Reidmeister-Schreier method,
followed by Tietze Transformations to simplify the presentation. This gave
$L=\scl{x,y:y^3=xyxy^{-1}xyx^{-1}yx=1}$ directly.

The composition series of $W$ contains the Janko group $J_1$ and the
projective special linear group $L_2(19)$~-- these are the quotients of $W$ by (respectively)
the normal subgroups $N'$ and $N''$ of Section~\ref{twoquos}.
The composition series cannot contain any other factors, since the order
of $W$ is exactly the product of the orders of $J_1$ and $L_2(19)$.
It follows that $N'\cong L_2(19)$ and $N''\cong J_1$.

Since $|W|=|N'|.|N''|$, and since $N'$ and $N''$ are both normal
in $W$, it follows that if $W=N'N''$ then $W$ is an internal direct product
$N'\times N''$. Note in any case that $N'N''$ is closed under the
group multiplication. To show $W=N'\times N''$
it therefore suffices to show that the generators of $W$ are found in $N'N''$.

Now $N'$ is generated by all conjugates of $\nu^3$, and
$N''$ by all conjugates of $\omega$. However,
$\nu^3=(\omega s_3)^3=s_3(s_3\omega s_3)\omega(s_3\omega s_3)$. Therefore,
the generator $s_3$ of $W$ may be expressed as a product
$\nu^3(s_3\omega s_3)\omega(s_3\omega s_3)$ of an element of $N'$ with
an element of $N''$. Since $s_2=(s_3s_2s_3s_2)s_3(s_2s_3s_2s_3)$, $s_1=(s_2s_1)s_2(s_1s_2)$
and $s_0=(s_1s_0s_1s_0)s_1(s_0s_1s_0s_1)$, that is, all generators of
$W$ are mutually conjugate, it follows that each generator of $W$ may
be written $\alpha'\alpha''$ for some conjugate $\alpha'$ of $\nu^3$ and
some conjugate $\alpha''$ of $(s_3\omega s_3)\omega(s_3\omega s_3)$.
Therefore $W=N'\times N''$ as required.
\ndpf

The twenty right
cosets of $L$ are $L$, $Ls_0$, $Ls_1=Ls_2$, $Ls_3$, $Ls_0s_1$, $Ls_0s_2$,
$Ls_0s_3$, $Ls_1s_3$, $Ls_3s_2$, $Ls_0s_1s_0$, $Ls_0s_1s_2$,
$Ls_0s_1s_3$, $Ls_0s_2s_1$, $Ls_0s_2s_3$, $Ls_0s_3s_2$, $Ls_1s_3s_2$,
$Ls_3s_2s_1$, $Ls_3s_2s_3$, $Ls_0s_1s_0s_2$ and $Ls_0s_3s_2s_3$.
Numbering the cosets in the order given yields a representation of the
permutation action of $W$ on these cosets, that is, a homomorphism from
$W$ to $S_{20}$:
$$s_0\mapsto(1,2) (3,6) (4,7) (5,10) (8,14) (9,15) (11,19) (12,16) (13,17) (18,20),$$
$$s_1\mapsto(1,3) (2,5) (4,8) (6,13) (7,12) (9,17) (10,15) (11,18) (14,19) (16,20),$$
$$s_2\mapsto(1,3) (2,6) (4,9) (5,11) (7,15) (8,16) (10,19) (12,14) (13,18) (17,20),$$
$$s_3\mapsto(1,4) (2,7) (3,8) (5,12) (6,14) (9,18) (10,16) (11,17) (13,19) (15,20).$$
This is not a faithful action of $W$ of course.
In fact, these permutations generate
a group isomorphic to $L_2(19)\cong W/N''\cong N'$.

Since $W\cong J_1\times L_2(19)$, it has a permutation presentation on
286 points, which may be derived from the presentations for $W'$ and $W''$ of
the previous section, via $s_i=s'_is''_i$. This presentation facilitates machine
computation. As before, see the notes at the end of Section~\ref{allquos}
for infomation on obtaining this presentation.

Note that $\nu^6=1$, since $\nu^6=(\nu^3)^2=(\omega(s_3\omega s_3))^3$ and
is therefore an element of $N'\cap N''$.

\section{The Remaining Quotients of $\CP$.}
\label{allquos}
\seq

By Theorem 2.7 of \cite{har6}, the semisparse subgroups $N$ of $W$ may be characterised by the property that
$N^w\cap H_0H_3\subseteq\{1,\omega\}$ for all $w\in W$. This is because any subgroup
satisfying this property is semisparse, and conversely
any semisparse subgroup of $W$ satisfies this property (see Theorem 2.7 of \cite{har6}).

The group $J_1\times L_2(19)$ has 1262 conjugacy classes of subgroups. It is relatively
straightforward to check them one by one to see if they satisfy the above
property. This was done, using Magma \cite{magma}. The program
took approximately 40 days of computing time, on two Intel Xeon processors running at
2GHz with 3Gb of RAM. The authors would like to suggest that
the results could have been obtained faster if they had troubled to optimise
the code better. In particular, in \cite{har6} it was noted that
if a C-group (such as $W$) satisfies the conditions of Theorem 2.7 of \cite{har6},
then whenever $N$ is semisparse in $W$, all its subgroups are also be semisparse.
This property was not used, but could have saved a significant amount of
computation time.

A total of 145 conjugacy classes of semisparse subgroups were
discovered, yielding 145 polytopes, most of them new. Further analysis of
the semisparse subgroups, to identify for example the facet types and
automorphism groups of the polytopes, was done using GAP version 4 release 3
(\cite{gap4}).

As just mentioned, if $N$ is semisparse
and $K\leq N$, then $K$ is also semisparse. Furthermore, $\CP/K$ is
a cover for $\CP/N$. For these reasons, it is important to know the
subgroup relations between the 145 semisparse subgroups of $W$. In particular,
it is useful to know the ``maximal'' semisparse subgroups, that is, those
semisparse subgroups which are not proper subgroups of any other.
These subgroups are important because the semisparse subgroups of $W$ are
exactly the subgroups of these maximal ones.

$W$ has 30 maximal semisparse subgroups.
These, with their generating sets, are listed
in Table \ref{maxgeo}. In the table, $v_1=\nu=(s_0s_1s_2)^5s_3$, $v_2=v_1s_2$, and
$v_6=v_5s_0=v_4s_1s_0=v_3s_0s_1s_0=v_2(s_1s_0)^2$. The notation used in the
``Group'' column is as follows: $G\times H$ is the direct product of $G$ with $H$,
and $G:H$ is a semidirect product. $G^k$ means the direct product of
$k$ copies of $G$. Also, as usual,
$C_k$ (or just `$k$') means the cyclic group of order $k$, $D_{2k}$ the dihedral
group of order $2k$, and $A_k$ and $S_k$ the alternating and symmetric
groups on $k$ points. Finally, as before,
$J_1$ is the first Janko group, and
$L_2(p)$ is the projective special linear group of rank 2 over $GF(p)$.
The numbering of the groups is in accordance with the ordering of the semisparse subgroups
as returned by the algorithms used by the authors.
This ordering is descending in the size of the group.

\begin{table}[ht]
\begin{center}
\begin{tabular}{|c|c|c|}
\hline
 No. & Group & Generators \\
 \hline
1 & $J_1$ &  $v_1^2,  (v_2 v_1)^2,  v_3^5,   (v_5 v_2 v_1)^2  $ \\
2 & $L_2(19)$ &  $(v_2 v_1^2)^2,  (v_6^2 v_4 v_1)^2  $ \\
3 & $19^2:9$ &  $v_1^2 v_2 v_4^2 v_6 v_3,  v_3^2 v_1 v_4 v_5^2 v_6 v_1  $ \\
4 & $19^2:9$ &  $v_5 v_1 v_5 v_4 v_1 v_2 v_4 v_1,  v_2^2 v_1 v_5 v_3^2 v_2 v_4$ \\
5 & $D_{38}\times 19:3$ &  $v_1 v_3 v_4^2 v_5 v_2^2 v_1,  v_1 v_6^2 v_2 v_3 v_5 v_2 v_4$ \\
6 & $19\times 11:10$ &  $v_1 v_3 v_6 v_5^2 v_3,  v_2^2 v_6^2 v_1 v_6 v_2$ \\
7 & $2^3\times 19:9$ &  $v_3 v_4 v_1 v_6 v_5 v_3 v_1,  v_5 v_6 v_5 v_1 v_2 v_5^2,   v_6 v_2^2 v_6 v_4 v_1 v_4 v_2$ \\
8 & $19\times A_5$ &  $v_5 v_6^2 v_3 v_5 v_1,  v_3^2 v_2 v_3 v_1 v_2 v_3,   v_3 v_5 v_2 v_6 v_5 v_3 v_6$ \\
9 & $19\times A_5$ &  $v_5 v_6^2 v_3 v_5 v_1,  v_2 v_3^3 v_2 v_3 v_1  $ \\
11 & $19\times 19:3$ &  $v_2 v_5 v_2 v_1 v_3 v_4 v_3,  v_4^2 v_1 v_6 v_5 v_3 v_2 v_1^2,   v_1 v_3 v_4 v_5 v_6 v_1 v_4 v_3 v_1  $ \\
13 & $5\times 19:9$ &  $v_1 v_3^6 v_1^2,  v_5 v_3 v_5 v_6 v_4 v_6^2 v_2 v_1  $ \\
16 & $2^2:(19:9)$ &  $v_4 v_6^2 v_2 v_3 v_4^2 v_2,  v_1 v_5 v_6 v_5^2 v_6 v_3 v_4$ \\
17 & $2^2:(19:9)$ &  $v_1 v_2 v_3 v_6 v_1 v_4 v_2 v_1,  v_4 v_6^2 v_2 v_3 v_4^2 v_2$ \\
19 & $19\times D_{30}$ &  $v_6 v_3 v_4 v_6 v_5^2,  v_6 v_2^3 v_6^2 v_1^3,   v_5 v_6^2 v_5 v_3 v_5 v_4 v_2 v_1  $ \\
20 & $C_{95}\times S_3$ &  $v_3^2 v_5^2 v_6 v_5 v_4 v_1,  (v_6 v_1)^2 v_2 v_5 v_4 v_5 v_1  $ \\
22 & $7\times A_5$ &  $v_1 v_5 v_4 v_2,  v_5 v_1 v_3 v_5 v_1 v_4 v_2^2 v_6$ \\
23 & $7\times A_5$ &  $v_1^3,  v_4 v_3 v_4 v_2 v_6 v_3 v_2^2 v_4,   v_4 v_6 v_5 v_4^2 v_3 v_5 v_3 v_4$ \\
25 & $19\times D_{20}$ &  $v_3 v_1^2 v_5^2 v_6 v_4 v_5,  v_5 v_6^2 v_3 v_5 v_4 v_2 v_5 v_4$ \\
29 & $(2^3:7)\times 5$ &  $(v_3 v_1)^2 v_5 v_2 v_1 v_4,  v_6^2 v_3 v_4 v_5^2 v_4 v_5$ \\
30 & $19\times A_4$ &  $v_3^6 v_1^3,  v_6 v_3^2 v_4 v_3^2 v_1 v_4 v_1  $ \\
33 & $11\times D_{20}$ &  $v_4 v_5 v_4 v_3 v_1^2,  v_2 v_1^3 v_2^4$ \\
42 & $(2^3:7)\times3$ &  $v_1 v_5 v_3 v_1 v_4 v_6^2 v_1^2,  v_3 v_6 v_2 v_6 v_4 v_6^3 v_1  $ \\
48 & $19\times S_3$ &  $v_6 v_4 v_3 v_4 v_1 v_5 v_2 v_1 v_6 v_1,   v_3 v_6 v_5 v_2 v_3 v_4^2 v_5 v_6 v_1  $ \\
54 & $5\times D_{22}$ &  $v_3 v_4 v_1 v_5 v_2^2 v_4,  v_1^2 (v_2^2 v_1)^2  $ \\
58 & $5\times D_{14}$ &  $v_6 v_1^2 v_2 v_5 v_6 v_3 v_4 v_1,  v_6 v_1^2 v_3 v_5^2 v_1^2 v_4 v_1  $ \\
64 & $7:9$ &  $v_5 v_3 v_2 v_5 v_6 v_1 (v_3 v_1)^2,   v_2^2 v_4 v_1^2 v_4 v_3 v_6 v_3 v_2$ \\
70 & $3\times D_{20}$ &  $v_2^3 v_6^2 v_1 v_6,  v_5 v_3 v_5 v_4 v_2 v_3 v_4 v_1 v_4 v_2,   v_1 v_5 v_1 v_3 v_5 v_2 v_1 v_6 v_5 v_6$ \\
81 & $3\times D_{14}$ &  $v_5 v_3 v_1 v_2 v_4 v_2 v_1 v_6 v_1,  (v_4 v_5 v_4)^3$ \\
95 & $5\times S_3$ &  $v_1 v_2 v_1 v_4 v_5^2 v_4 v_6 v_4,   v_5 v_3 v_5 v_4 v_5 v_6 v_4 v_2 v_6 v_1^2  $ \\
127 & $C_{10}$ &  $v_1 v_2 v_3 v_1 v_2 v_4 v_1^2 v_3^2 v_1  $  \\
 \hline
 \end{tabular}
 \end{center}
 \caption{The Maximal Semisparse Subgroups of $W=[\{5,3\},\{3,5\}_5]$.}
 \label{maxgeo}
 \end{table}

Of the 145 semisparse subgroups, 70 yield section regular polytope. Of these,
69 have dodecahedral facets, and one (namely, number 1) is Coxeter's 57-cell,
self-dual, with hemidodecahedral facets.
A list of these may be found in Tables \ref{secreg1} and \ref{secreg2}. In
those tables is listed the
number of the group, the isomorphism type of the group, its maximal proper
semisparse subgroups in $W$, the number of facets of the polytope, and the
isomorphism type of its automorphism group. Note that group number 40 is not isomorphic
to either group 39 or 38, nor any other group labelled $19:9$ in this paper.
It is perhaps better described as a non-split extension of a normal subgroup
$C_{57}$ by a group of order $3$.

\begin{table}
\begin{center}
\begin{tabular}{|c|l|l|c|l|}
\hline
No. & Group & Subgroups & \#Facets & $\Aut(\CP)$ \\
\hline
1 & $J_1$ &  18, 41, 44, 45, 51, 65, 80 & 57 & $L_2(19)$ \\
2 & $L_2(19)$ &  38, 66, 67, 103, 109 & 1463 & $J_1$ \\
3 & $19^2:9$ &  10, 39, 40 & 1540 & $C_6$ \\
4 & $19^2:9$ &  10, 39, 40 & 1540 & $C_6$ \\
10 & $19\times 19:3$ &  26, 71, 72 & 4620 & $3^2\times 2$ \\
11 & $19\times 19:3$ &  26, 73, 74 & 4620 & $C_{18}$ \\
12 & $19\times 11:5$ &  34, 56, 76 & 4788 & $C_{18}$ \\
13 & $5\times 19:9$ &  38, 56, 78 & 5852 & $D_{12}$ \\
22 & $7\times A_5$ &  57, 62, 66, 82 & 11913 & $C_6$ \\
23 & $7\times A_5$ &  57, 61, 67, 82 & 11913 & $C_6$ \\
26 & $19^2$ &  106, 107, 108 & 13860 & $C_6\times C_9$ \\
28 & $C_{285}$ &  56, 74, 112 & 17556 & $2^2\times C_9$ \\
30 & $19\times A_4$ &  58, 72, 116 & 21945 & $C_6$ \\
33 & $11\times D_{20}$ &  52, 53, 55, 79, 103 & 22743 & $C_{10}$ \\
34 & $C_{209}$ &  107, 119 & 23940 & $C_{90}$ \\
38 & $19:9$ &  71, 128 & 29260 & $J_1$ \\
39 & $19:9$ &  71, 129 & 29260 & $3\times D_{10}$ \\
40 & $19:9~(\cong 57.3)$ &  72, 129 & 29260 & $C_6$ \\
48 & $19\times S_3$ &  72, 86, 134 & 43890 & $C_6$ \\
52 & $11\times D_{10}$ &  77, 100, 121 & 45486 & $2^2\times 5$ \\
53 & $11\times D_{10}$ &  77, 100, 122 & 45486 & $2^2\times 5$ \\
55 & $C_{110}$ &  77, 100, 126 & 45486 & $2^2\times 5$ \\
56 & $C_{95}$ &  107, 137 & 52668 & $9\times D_{12}$ \\
57 & $7\times A_4$ &  96, 102, 116 & 59565 & $C_6$ \\
58 & $2^2\times 19$ &  86, 139 & 65835 & $3^2\times 2$ \\
61 & $7\times D_{10}$ &  89, 114, 121 & 71478 & $2^2\times 3$ \\
62 & $7\times D_{10}$ &  89, 114, 122 & 71478 & $2^2\times 3$ \\
64 & $7:9$ &  102, 129 & 79420 & $C_6$ \\
66 & $A_5$ &  116, 122, 134 & 83391 & $J_1$ \\
67 & $A_5$ &  116, 121, 134 & 83391 & $J_1$ \\
71 & $19:3$ &  107, 141 & 87780 & $3\times J_1$ \\
72 & $C_{57}$ &  106, 141 & 87780 & $C_6\times S_3$ \\
73 & $19:3$ &  106, 142 & 87780 & $2\times L_2(19)$ \\
74 & $C_{57}$ &  107, 142 & 87780 & $9\times D_{20}$ \\
76 & $11:5$ &  119, 137 & 90972 & $2\times L_2(19)$ \\
77 & $C_{55}$ &  119, 138 & 90972 & $2^3\times 5$ \\
78 & $C_{45}$ &  111, 128 & 111188 & $2\times D_{12}$ \\
79 & $2^2\times 11$ &  100, 139 & 113715 & $C_{30}$ \\
82 & $7\times S_3$ &  102, 114, 134 & 119130 & $C_6$ \\
86 & $C_{38}$ &  106, 144 & 131670 & $C_6\times D_{10}$ \\
\hline
\end{tabular}
\caption{Semisparse Subgroups Corresponding to Section Regular Polytopes - Part 1}
\label{secreg1}
\end{center}
\end{table}

\begin{table}
\begin{center}
\begin{tabular}{|c|l|l|c|l|}
\hline
No. & Group & Subgroups & \#Facets & $\Aut(\CP)$ \\
\hline
89 & $C_{35}$ &  131, 138 & 142956 & $2^3\times 3$ \\
96 & $2^2\times 7$ &  114, 139 & 178695 & $3^2\times 2$ \\
100 & $C_{22}$ &  119, 144 & 227430 & $C_{10}\times D_{10}$ \\
101 & $7:3$ &  131, 142 & 238260 & $2\times L_2(19)$ \\
102 & $C_{21}$ &  131, 141 & 238260 & $C_6\times S_3$ \\
103 & $D_{20}$ &  121, 122, 126, 139 & 250173 & $J_1$ \\
106 & $19$ &  145 & 263340 & $C_6\times L_2(19)$ \\
107 & $19$ &  145 & 263340 & $C_9\times J_1$ \\
108 & $19$ &  145 & 263340 & $19:3$ \\
109 & $9:2$ &  128, 134 & 277970 & $J_1$ \\
111 & $C_{15}$ &  137, 141 & 333564 & $S_3\times D_{12}$ \\
112 & $C_{15}$ &  137, 142 & 333564 & $2^2\times L_2(19)$ \\
114 & $C_{14}$ &  131, 144 & 357390 & $C_6\times D_{10}$ \\
116 & $A_4$ &  139, 141 & 416955 & $J_1$ \\
119 & $11$ &  145 & 454860 & $C_{10}\times L_2(19)$ \\
121 & $D_{10}$ &  138, 144 & 500346 & $2\times J_1$ \\
122 & $D_{10}$ &  138, 144 & 500346 & $2\times J_1$ \\
126 & $C_{10}$ &  138, 144 & 500346 & $2\times J_1$ \\
127 & $C_{10}$ &  137, 144 & 500346 & $D_{10}\times D_{12}$ \\
128 & $C_9$ &  141 & 555940 & $2\times J_1$ \\
129 & $C_9$ &  141 & 555940 & $S_3\times D_{10}$ \\
131 & $7$ &  145 & 714780 & $C_6\times L_2(19)$ \\
134 & $S_3$ &  141, 144 & 833910 & $J_1$ \\
137 & $5$ &  145 & 1000692 & $D_{12}\times L_2(19)$ \\
138 & $5$ &  145 & 1000692 & $2^2\times J_1$ \\
139 & $2^2$ &  144 & 1250865 & $3\times J_1$ \\
141 & $3$ &  145 & 1667820 & $S_3\times J_1$ \\
142 & $3$ &  145 & 1667820 & $D_{20}\times L_2(19)$ \\
144 & $2$ &  145 & 2501730 & $D_{10}\times J_1$ \\
145 & $1$ & -- & 5003460 & $J_1\times L_2(19)$ \\
\hline
\end{tabular}
\caption{Semisparse Subgroups Corresponding to Section Regular Polytopes - Part 2}
\label{secreg2}
\end{center}
\end{table}

The remaining 75 groups (or polytopes)
are listed in Tables \ref{remain1} and \ref{remain2}.
Those two tables list similar information to that given for the section regular
polytopes, the key difference being in the information given about the
facets of the polytopes. An entry in the ``Facets'' columns of the form $D^dH^h$
means that the polytope has, as facets, $d$ dodecahedra and $h$ hemidodecahedra.

It should not seem surprising or contradictory that (for example) factoring
out by a single element of order $2$ (namely $\omega$) yields a polytope
with $3420$ hemi-dodecahedral facets (polytope number 143). For example,
factoring the cube, with group $\scl{\sigma_0,\sigma_1,\sigma_2}$, by the semisparse
subgroup $\{1,(\sigma_0\sigma_1)^2\}$ yields a digonal prism, a polytope with
not one, but two digons corresponding to ``opposite'' squares of the original cube.
Likewise here, if a ``cusp'' is introduced into the ``space'' occupied by $\CP$
by factoring out $\omega$, this forces another $3419$ cusps to also form.
The elements of $N'=\scl{\scl{\nu^3}}\cong L_2(19)$ permute among the cusps.
The automorphism group $A_5$ of the hemi-icosahedron also acts on this ``space'' by
rotating it around the cusps, just as the group of the digon acts on the
digonal prism in a way that maps each digon to itself. This is why the automorphism
group of polytope 143 is $A_5\times L_2(19)$.

\begin{table}
\begin{center}
\begin{tabular}{|c|l|l|c|l|}
\hline
No. & Group & Subgroups & Facets & $\Aut(\CP)$ \\
\hline
5 & $19:3\times D_{38}$ &  10, 14, 46, 47 & $D^{2280}H^{60}$ & $3^2$ \\
6 & $19\times 11:10$ &  12, 24, 35, 51 & $D^{2376}H^{36}$ & $C_9$ \\
7 & $2^3\times 19:9$ &  15, 21, 60 & $D^{3640}H^{35}$ & $7:3$ \\
8 & $19\times A_5$ &  31, 36, 49, 68 & $D^{4344}H^{90}$ & $C_{18}$ \\
9 & $19\times A_5$ &  31, 36, 50, 69 & $D^{4344}H^{90}$ & $C_9$ \\
14 & $19\times D_{38}$ &  26, 84, 85 & $D^{6840}H^{180}$ & $3\times C_9$ \\
15 & $(2^2)\times (19:9)$ &  27, 32, 87 & $D^{7300}H^{30}$ & $C_6$ \\
16 & $(2^2):(19:9)$ &  32, 39, 88 & $D^{7300}H^{30}$ & $C_6$ \\
17 & $(2^2):(19:9)$ &  32, 39, 88 & $D^{7300}H^{30}$ & $C_6$ \\
18 & $L_2(11)$ &  68, 69, 76, 115 & $D^{7296}H^{570}$ & $L_2(19)$ \\
19 & $19 \times D_{30}$ &  28, 36, 50, 92 & $D^{8688}H^{180}$ & $C_{18}$ \\
20 & $C_{95}\times S_3$ &  28, 35, 49, 93 & $D^{8760}H^{36}$ & $C_{18}$ \\
21 & $2^2\times (19:3)$ &  32, 43, 98 & $D^{10920}H^{105}$ & $3\times 7:3$ \\
24 & $19\times D_{22}$ &  34, 85, 99 & $D^{11880}H^{180}$ & $C_{45}$ \\
25 & $19\times D_{20}$ &  35, 36, 37, 59, 104 & $D^{13068}H^{198}$ & $C_9$ \\
27 & $2\times 19:3$ &  38, 47, 110 & $D^{14620}H^{20}$ & $A_5$ \\
29 & $5\times (2^3:7)$ &  75, 83, 89 & $D^{17784}H^{171}$ & $2^2\times 3$ \\
31 & $19\times A_4$ &  59, 74, 117 & $D^{21900}H^{90}$ & $C_{18}$ \\
32 & $2^2\times (19:3)$ &  47, 59, 118 & $D^{21900}H^{90}$ & $3^2\times 2$ \\
35 & $C_{190}$ &  56, 85, 120 & $D^{26316}H^{36}$ & $C_{18}$ \\
36 & $19\times D_{10}$ &  56, 85, 123 & $D^{26244}H^{180}$ & $C_{18}$ \\
37 & $19\times D_{10}$ &  56, 85, 125 & $D^{26244}H^{180}$ & $S_3\times C_9$ \\
41 & $8:(7:3)$ &  75, 97, 101 & $D^{29640}H^{285}$ & $L_2(19)$ \\
42 & $(2^3:7)\times 3$ &  75, 98, 102 & $D^{29640}H^{285}$ & $3\times S_3$ \\
43 & $19\times 2^3$ &  59, 130 & $D^{32760}H^{315}$ & $9\times 7:3$ \\
44 & $2\times A_5$ &  68, 97, 104, 115 & $D^{40812}H^{1767}$ & $L_2(19)$ \\
45 & $19:6$ &  73, 84, 132 & $D^{43320}H^{1140}$ & $L_2(19)$ \\
46 & $3\times D_{38}$ &  72, 84, 133 & $D^{43320}H^{1140}$ & $3\times S_3$ \\
47 & $2\times 19:3$ &  71, 85, 133 & $D^{43860}H^{60}$ & $3\times A_5$ \\
49 & $19\times S_3$ &  74, 85, 135 & $D^{43800}H^{180}$ & $9\times D_{10}$ \\
50 & $19\times S_3$ &  74, 85, 136 & $D^{43800}H^{180}$ & $C_{18}$ \\
51 & $11:10$ &  76, 99, 120 & $D^{45144}H^{684}$ & $L_2(19)$ \\
54 & $5\times D_{22}$ &  77, 99, 124 & $D^{45144}H^{684}$ & $2^2\times 5$ \\
59 & $2^2\times 19$ &  85, 140 & $D^{65700}H^{270}$ & $C_6\times C_9$ \\
60 & $2^3\times 9$ &  87, 98 & $D^{69160}H^{665}$ & $2\times 7:3$ \\
63 & $5\times D_{14}$ &  89, 113, 124 & $D^{71136}H^{684}$ & $2^2\times 3$ \\
65 & $S_3\times D_{10}$ &  92, 93, 94, 104, 115 & $D^{82080}H^{2622}$ & $L_2(19)$ \\
68 & $A_5$ &  117, 123, 135 & $D^{82536}H^{1710}$ & $2\times L_2(19)$ \\
69 & $A_5$ &  117, 123, 136 & $D^{82536}H^{1710}$ & $L_2(19)$ \\
70 & $3\times D_{20}$ &  90, 91, 104, 118 & $D^{82764}H^{1254}$ & $S_3$ \\
\hline
\end{tabular}
\caption{The Remaining Semisparse Subgroups - Part 1}
\label{remain1}
\end{center}
\end{table}

\begin{table}
\begin{center}
\begin{tabular}{|c|l|l|c|l|}
\hline
No. & Group & Subgroups & Facets & $\Aut(\CP)$ \\
\hline
75 & $(2^3):7$ &  130, 131 & $D^{88920}H^{855}$ & $3\times L_2(19)$ \\
80 & $7:6$ &  101, 113, 132 & $D^{118560}H^{1140}$ & $L_2(19)$ \\
81 & $3\times D_{14}$ &  102, 113, 133 & $D^{118560}H^{1140}$ & $3\times S_3$ \\
83 & $2^3\times 5$ &  105, 130 & $D^{124488}H^{1197}$ & $2^2\times 7:3$ \\
84 & $D_{38}$ &  106, 143 & $D^{129960}H^{3420}$ & $3\times L_2(19)$ \\
85 & $C_{38}$ &  107, 143 & $D^{131580}H^{180}$ & $9\times A_5$ \\
87 & $2^2\times 9$ &  110, 118 & $D^{138700}H^{570}$ & $2^2\times 3$ \\
88 & $2^2:9$ &  118, 129 & $D^{138700}H^{570}$ & $C_6$ \\
90 & $C_{30}$ &  111, 120, 133 & $D^{166668}H^{228}$ & $D_{12}$ \\
91 & $3\times D_{10}$ &  111, 123, 133 & $D^{166212}H^{1140}$ & $D_{12}$ \\
92 & $D_{30}$ &  112, 123, 136 & $D^{165072}H^{3420}$ & $2\times L_2(19)$ \\
93 & $5\times S_3$ &  112, 120, 135 & $D^{166440}H^{684}$ & $2\times L_2(19)$ \\
94 & $3\times D_{10}$ &  112, 125, 132 & $D^{166212}H^{1140}$ & $2\times L_2(19)$ \\
95 & $5\times S_3$ &  124, 135 & $D^{166440}H^{684}$ & $2^2\times D_{10}$ \\
97 & $2\times A_4$ &  117, 130, 132 & $D^{207480}H^{1995}$ & $L_2(19)$ \\
98 & $2^3\times 3$ &  118, 130 & $D^{207480}H^{1995}$ & $S_3\times 7:3$ \\
99 & $D_{22}$ &  119, 143 & $D^{225720}H^{3420}$ & $5\times L_2(19)$ \\
104 & $D_{20}$ &  120, 123, 125, 140 & $D^{248292}H^{3762}$ & $L_2(19)$ \\
105 & $2^2\times 5$ &  124, 140 & $D^{249660}H^{1026}$ & $2^3\times 3$ \\
110 & $C_{18}$ &  128, 133 & $D^{277780}H^{380}$ & $2\times A_5$ \\
113 & $D_{14}$ &  131, 143 & $D^{355680}H^{3420}$ & $3\times L_2(19)$ \\
115 & $D_{12}$ &  132, 135, 136, 140 & $D^{414960}H^{3990}$ & $L_2(19)$ \\
117 & $A_4$ &  140, 142 & $D^{416100}H^{1710}$ & $2\times L_2(19)$ \\
118 & $2^2\times 3$ &  133, 140 & $D^{416100}H^{1710}$ & $C_6\times S_3$ \\
120 & $C_{10}$ &  137, 143 & $D^{500004}H^{684}$ & $2\times L_2(19)$ \\
123 & $D_{10}$ &  137, 143 & $D^{498636}H^{3420}$ & $2\times L_2(19)$ \\
124 & $C_{10}$ &  138, 143 & $D^{500004}H^{684}$ & $2^2\times A_5$ \\
125 & $D_{10}$ &  137, 143 & $D^{498636}H^{3420}$ & $S_3\times L_2(19)$ \\
130 & $2^3$ &  140 & $D^{622440}H^{5985}$ & $(7:3)\times L_2(19)$ \\
132 & $C_6$ &  142, 143 & $D^{833340}H^{1140}$ & $2\times L_2(19)$ \\
133 & $C_6$ &  141, 143 & $D^{833340}H^{1140}$ & $S_3\times A_5$ \\
135 & $S_3$ &  142, 143 & $D^{832200}H^{3420}$ & $D_{10}\times L_2(19)$ \\
136 & $S_3$ &  142, 143 & $D^{832200}H^{3420}$ & $2\times L_2(19)$ \\
140 & $2^2$ &  143 & $D^{1248300}H^{5130}$ & $C_6\times L_2(19)$ \\
143 & $2$ &  145 & $D^{2500020}H^{3420}$ & $A_5\times L_2(19)$ \\
\hline
\end{tabular}
\caption{The Remaining Semisparse Subgroups - Part 2}
\label{remain2}
\end{center}
\end{table}

As mentioned earlier, some auxiliary information for this article is available.
The auxiliary information is in the form of a file
({\tt g.txt}) containing GAP commands that construct the permutation presentations
for the groups of the polytopes $\CP'$ and $\CP''$, and combine them into
a presentation (on 286 points) for $W$.

The file also contains permutation presentations for representatives of
all the 1262 conjugacy classes of subgroups of $W\cong J_1\times L_2(19)$.
The file defines a list ({\tt cl}) of indices to identify which of these groups
are semisparse in $W$. Importing this file into GAP (via, for example,
a {\tt Read} command) places in GAP's workspace a list ({\tt geo})
consisting of one representative of each conjugacy class of semisparse subgroups
of $W$.

The file {\tt g.txt} is available via a link from \cite{535aux}. If that
web page should move to a different URL, it should in any case
be locatable by searching the web for the title of this article.

The first author would like to acknowledge and thank Jesus Christ, through whom
all things were made, for the encouragement, inspiration, and occasional hint
that were necessary to complete this article. The second author, however,
specifically disclaims this acknowledgement.

The second author would like to acknowledge financial support from the
Belgian National Fund for Scientific Research which made this project
feasible.

Both authors acknowledge with gratitude the helpful comments made by
Buekenhout on a draft of this paper.

\end{document}